\documentclass[reqno]{amsart}

\usepackage{color}

\usepackage{amsthm}
\usepackage{amssymb}
\usepackage{amsxtra}
\usepackage{graphics}
\usepackage{graphicx}

\title[A Lax-Wendroff type theorem]{A Lax-Wendroff type theorem for unstructured quasiuniform grids}
\author[V. Elling]{Volker Elling}
\address{Volker Elling\\SCCM Program\\Gates Building 2B\\Stanford University\\Stanford, CA 94305-9025}
\email{velling@stanford.edu}
\urladdr{http://www.stanford.edu/\textasciitilde velling}
\thanks{This material is based upon work supported by an SAP/Stanford Graduate Fellowship and by the National Science Foundation under
Grant no. DMS 0104019. Any opinions, findings, and conclusions or recommendations expressed in this material are those of the author and
do not necessarily reflect the views of the National Science Foundation.}


\theoremstyle{plain}
\newtheorem{thm}{Theorem}
\newtheorem{lemma}{Lemma}
\newtheorem{cor}{Corollary}
\theoremstyle{definition}

\newtheorem{exa}{Example}
\theoremstyle{remark}

\newcommand{\aint}{\hbox to-1.5mm{--}\int}
\newcommand{\defm}[1]{{\it #1}}
\newcommand{\loc}{\operatorname{loc}}
\newcommand{\stn}{\operatorname{stn}}
\newcommand{\supp}{\operatorname{supp}}
\newcommand{\interior}{\operatorname{int}}
\newcommand{\sys}[1]{\mathcal{#1}}
\newcommand{\R}{\mathbb{R}}
\newcommand{\RC}{\R^{d+1}_+}
\newcommand{\Rd}{\R^d}
\newcommand{\Z}{\mathbb{Z}}
\newcommand{\ZC}{\Nnull\times\Z^d}

\newcommand{\Nnull}{\mathbb{N}_0}

\newcommand{\diam}{\operatorname{diam}}
\newcommand{\abb}[3]{#1:#2\rightarrow #3}
\newcommand{\ef}[1]{\frac{1}{#1}}
\newcommand{\Div}{\operatorname{div}}
\newcommand{\rarrow}{\bgroup\rightarrow\egroup}
\newcommand{\Nd}{\{0,\dotsc,d\}}
\newcommand{\Eta}{E}

\newcommand{\subeq}[2]{\mathord{\underbrace{\mathop{#1}}_{#2}}}

\newcommand{\topanno}[2]{\mathrel{\mathop{#2}^{#1}}}

\renewcommand{\Vec}[1]{#1}

\newcommand{\ucon}[1]{YUCK}

\makeatletter
\newcounter{aconst}

\makeatother

\def\doprivate{01}

\begin{document}

\begin{abstract}
    \parindent=0mm%
    \parskip=3mm%
    A well-known theorem of Lax and Wendroff states that if 
    the sequence of approximate solutions to a system of hyperbolic conservation laws generated by a conservative consistent numerical scheme 
    converges boundedly a.e.\ as the mesh parameter goes to zero, then the limit is a weak solution of the system. 
    Moreover, if the scheme satisfies a discrete entropy inequality as well, the limit is an entropy solution. 
    The original theorem applies to uniform Cartesian grids; 
    this article presents a generalization for quasiuniform grids (with Lipschitz-boundary cells)
    uniformly continuous inhomogeneous numerical fluxes and nonlinear inhomogeneous sources. The added generality allows a discussion of novel 
    applications like local time stepping, grids with moving vertices and conservative remapping. 
    A counterexample demonstrates that the theorem is not valid for arbitrary non-quasiuniform grids.
\end{abstract}

\maketitle


\parindent=0mm%
\parskip=3mm%

\section{Introduction}

Consider the Cauchy problem for systems of first-order conservation laws
\begin{alignat}{1}
    \sum_{i=0}^d\frac{d}{dy_i}f_i(u(y),y) &= p(u(y),y)\qquad(y\in\RC), \label{conservation-law} \\
    u(0,x) &= u_0(x)\qquad(x\in\R^d), \label{inicond}
\end{alignat}
where $\RC:=(0,\infty)\times\R^d$,
$\abb{u}{\RC}{P}$ (where $P\subset\R^m$ is a bounded open subset of the set of physically admissible values), 
$\Vec f=(f_0,\dotsc,f_d)'$ with $f_0(w)=w$ and smooth \defm{fluxes} $\abb{f_i=(f_{i1},\dotsc,f_{im})'}{P\times\overline{\RC}}{\R^m}$ ($i=1,\dotsc,d$),
smooth \defm{source} $p=(p_1,\dotsc,p_m)':P\times\overline{\RC}\rightarrow\R^m$ and \defm{initial values} $u_0:\R^d\rightarrow P$.

For the analysis of initial-value problems it is common to separate the time variable and the spatial variable(s); however, for the purposes of
the Lax-Wendroff theorem there is no benefit in distinguishing them. For brevity of notation we collect them in the single vector
$y=(t,x')'\in\RC$. $x$ will be used for coordinates in $\R^d$, $V$ resp.\ $S$ for the 
$(d+1)$-dimensional resp.\ $d$-dimensional Hausdorff measure.

It is well-known that even for smooth initial values $u_0$, there need not exist a smooth solution to \eqref{conservation-law}
for all $y_0>0$. It is necessary to extend the search to weak solutions, i.e.\ to $u\in L^1(\RC;P)$ that satisfy
\begin{alignat}{1}
    -\int_{\RC}\sum_{i=0}^df_i(u(y),y)\frac{\partial\phi}{\partial y_i}(y)dV(y) &= \int_{\R^d}\phi(0,\cdot)u_0~dS + \int_{\RC}p(u(y),y)\phi(y)dV(y) \label{weak-conservation-law}
\end{alignat}
for all $\phi\in C^\infty_c(\overline{\RC})$. Since there can be more than one weak solution, an entropy condition is needed
to select the ``physical'' one: let $\Vec\eta=(\eta_0,\dotsc,\eta_d)'$, $\eta_i:P\times\overline{\RC}\rightarrow\R$ smooth, $\eta_0$ (the \defm{entropy})
strictly convex and $\eta_1,\dotsc,\eta_d$ (\defm{entropy fluxes}) such that
\begin{equation}
    \frac{\partial\eta_i}{\partial u_\alpha} = \sum_{\beta=1}^m\frac{\partial\eta_0}{\partial u_\beta}\frac{\partial f_{i\beta}}{\partial u_\alpha}
    \qquad(i=1,\dotsc,d,~\alpha=1,\dotsc,m);
\end{equation}
$\Vec\eta$ is called \defm{entropy/entropy flux pair}.
The entropy condition is 
\begin{alignat}{1}
    & \sum_{i=0}^d\frac{d}{dy_i}\eta_i(u(y),y) \leq 
    \subeq{\sum_{\beta=1}^m\frac{\partial\eta_0}{\partial u_\beta}(u(y),y)p_\beta(u(y),y) + \sum_{i=0}^d\frac{\partial\eta_i}{\partial y_i}(u(y),y)}{=:g(u(y),y)}
    \label{entropy-condition}
\end{alignat}
which is meant to hold in the weak sense, i.e.\ for all \emph{nonnegative} $\phi\in C^\infty_c(\overline{\RC})$
\begin{alignat}{1}
    & -\int_{\RC}\sum_{i=0}^d\eta_i(u(y),y)\frac{\partial\phi}{\partial y_i}(y)dV(y)
    \leq \int_{\Rd}\eta_0(u_0(\cdot))\phi(0,\cdot)dV
    + \int_{\RC}\phi(y)g(u(y),y)dV.
    \label{weak-entropy-condition}
\end{alignat}
(Note: whether this entropy condition is sufficient to guarantee uniqueness is not known, except for some special cases.)

The rest of this section is limited to the case of conservation laws without sources; the presence of sources, as in reactive flow,
poses additional difficulties.

The classical proof that a sequence $(u^h)$ of numerical approximations to \eqref{conservation-law} converges to 
an entropy solution proceeds as follows: the properties of the numerical scheme (e.g.\ a monotone conservative scheme with consistent 
homogeneous (i.e.\ $y$-independent) fluxes on a uniform Cartesian grid, see \cite{harten-hyman-lax}, \cite{crandall-majda} and \cite{crandall-tartar}) 
guarantee that $(u^h)$ is bounded in $L^\infty$ and $TV$ (the space 
of functions with bounded variation in the sense of Tonelli-Cesari). This implies that some subsequences converge pointwise almost everywhere (in the presence of uniform $L^\infty$ boundedness equivalent to $L^1_{\loc}$ convergence);
the Lax-Wendroff theorem (see \cite{lax-wendroff}) proves that the limits of these subsequences are indeed weak 
solutions of \eqref{conservation-law}. Moreover, if the $u^h$ satisfy \defm{discrete entropy inequalities}, then the limits must be 
entropy solutions of \eqref{conservation-law}. Whenever entropy solutions are unique, the entire sequence $(u^h)$ must converge to the 
entropy solution. Positive uniqueness results are available in special cases (see \cite{kruzkov} for scalar conservation laws in multiple dimensions 
or \cite{bressan-lefloch} for 1D system entropy solutions with small total variation and some other restrictions), but see \cite{elling-nonuniqueness}
for a possible counterexample for 2D Euler system solutions.

In many cases, $L^1_{\loc}$ precompactness is difficult to prove --- and might be false ---, e.g.\ for unstructured grids or
higher-order schemes for scalar conservation laws, not to mention schemes for \emph{systems} of conservation laws. 
(For this reason, techniques based on measure-valued solutions which require $L^\infty$ boundedness, but not $L^1_{\loc}$ precompactness,
have been developed and successfully applied to the scalar case in \cite{coquel-lefloch-i}, \cite{coquel-lefloch-ii};
see also \cite{noelle-mv-irregular} for irregular grids.)
However, \cite{cockburn-coquel-lefloch} have generalized earlier work by \cite{kuznetsov} (see also \cite{sanders}) to a large class of unstructured 
grids.
They prove $L^1$ convergence of order $\ef{4}$ to the entropy solution, for monotone numerical 
fluxes with antidiffusive modifications; the modifications allow for higher order in regions where the entropy solution is smooth.
Although these Kuznetsov-type proofs yield convergence without resort to the Lax-Wendroff theorem, they demonstrate that the preconditions
of the Lax-Wendroff theorem are satisfied more often than previously thought.

More importantly, while \emph{rigorous} proofs of convergence are limited to special cases, 
\emph{observing} actual output of good numerical 
schemes suggests that bounded a.e.\ convergence is rather common, even for important systems like compressible gas dynamics. 
In this sense, the Lax-Wendroff theorem has important 
\emph{heuristic} value: it guarantees that 
the limit, if there is one, is a weak solution; moreover, in the presence of a discrete entropy condition, it guarantees that the limit is an entropy solution. 
\if\doprivate%
In fact, one could say, tongue-in-cheek, that one could even verify $TV$ boundedness, hence $L^1$ compactness, by inspection, in absence of obvious
$L^1$ convergence...
\fi
Finally, the Lax-Wendroff theorem serves as a theoretical motivation for focusing on conservative schemes with consistent fluxes 
(however, occasionally nonconservative schemes are used in practice).

While the Kuznetsov-type proof in \cite{cockburn-coquel-lefloch} applies to a large class of unstructured meshes, it relies 
strongly on special properties of scalar conservation laws; the same holds for techniques based on weak convergence and 
measure-valued solutions. It seems that only the Lax-Wendroff theorem provides at least a partial result for systems.

The original Lax-Wendroff theorem requires a 1D uniform Cartesian grid, continuous fluxes, 
$L^1_{\loc}$ precompactness and $L^\infty$ boundedness. 
LeVeque \cite{leveque-book} Section 12.4 simplifies the proof, at the cost of requiring locally Lipschitz-continuous numerical fluxes 
and $TV$ boundedness. 
\cite{kroener-rokyta-wierse} present a proof for 2D polygonal meshes, locally Lipschitz-continuous numerical fluxes, $L^\infty$ boundedness,
$L^1_{\loc}$ precompactness and an explicit CFL condition; however, their assumptions (2.3) and (2.4) about the mesh seem restrictive. 
More general triangular meshes are covered by Proposition 4.4.1 in \cite{godlewski-raviart}; 
it might be possible to extend their proof technique to polygonal meshes.
With straightforward modifications to statement and proof, all of these results and proofs apply to an arbitrary number of dimensions; however,
none of them seem to generalize into other directions easily. 
This article considers a quasiuniform mesh with no other geometric restrictions, (uniformly) continuous inhomogeneous numerical fluxes and nonlinear
inhomogeneous source terms.

Only the Cauchy problem is discussed; boundary conditions pose many open research problems, both 
theoretically and numerically. 
Even in benign cases where the flux is completely prescribed and independent of the solution near the boundary (as in supersonic inflow), 
one needs to make additional assumptions about the convergence of the numerical solution near the boundary which are not implied by mere 
boundedly a.e.\ convergence.

Section \ref{section:definitions} introduces the grids, numerical fluxes and numerical sources and the conditions imposed on them; 
this rather abstract framework is 
illustrated by a simple example in Section \ref{section:example}. Section \ref{section:fully-discrete} contains statement and proof of the 
generalized Lax-Wendroff theorem (Theorem \ref{thm:lax-wendroff-ich}).
Section \ref{section:counterexample} provides a counterexample that explains why Theorem \ref{thm:lax-wendroff-ich} does not always hold
for non-quasiuniform grids. The newfound generality enables theoretical discussion of some numerical techniques and applications 
in Section \ref{section:novel-applications}.


\section{Notation and assumptions}

\label{section:definitions}

\subsection{Landau symbols}
Two sequences of grids will be used: unstructured grids with parameter $h$, and uniform Cartesian grids with parameter $H$.
An expression $A$ is said to be $O(B)$ ($B$ some other expression) if there is some constant $c$, independent of
\[ C,N,F,h,H,\epsilon,\rho,k,w, \]
so that
\[ A \leq cB \]
as long as $h,H\in(0,1]$ and as long as
\begin{alignat}{1}
    \rho &:= \frac{h}{H} < \ef{2} \label{eq:rho-def}
\end{alignat}

$A$ is said to be $\Omega(B)$ if $B$ is $O(A)$. 

We say that an expression is $o_{\rho,\epsilon}(1)$ if, for any fixed values of $\rho,\epsilon>0$ (and $h:=\rho H$) it converges to $0$ as $H\downarrow 0$.

\subsection{Grids} 
\label{subsection:grids}
For any $h>0$, let $\sys C^h$ be a system of closed subsets (called \defm{cells}) of $\RC$ with pairwise disjoint interiors so that
\[ \bigcup_{C\in\sys C^h}C=\overline{\RC}. \]
We require the cells to have Lipschitz boundaries; this is more than weak enough for all conceivable numerical meshes.
For $C,N\in\sys C^h$ with $S(C\cap N)>0$, let $C\rarrow N$ denote the ordered pair $(C,N)$; depending on the context it will refer to $C\cap N$ instead.
The \defm{unit normal} $\Vec n_{C\rarrow N}(y)$ in each point $y\in C\cap N$ is fixed as pointing into $N$.
The $C\rarrow N$ are called \defm{interior faces}; the other class of faces consists of \defm{initial faces} $C\rarrow\partial$, $\partial\rarrow C$
(where $C\cap(\{0\}\times\R^d)\neq\emptyset$).
$\partial\rarrow C$ will sometimes refer to $C\cap(\{0\}\times\R^d)$, with unit normal $(1,0,\dotsc,0)'$; 
$C\rarrow\partial$ will refer to the same surface with opposite unit normal. Define $\hat{\sys C}^h:=\sys C^h\cup\{\partial\}$.
Let $\sys F^h$ be the set of interior faces, $\hat{\sys F}^h$ the set of all faces.

To ``define'' the mesh parameter $h$, require
\begin{equation}
    \diam C\leq h\qquad(C\in\sys C^h); \label{eq:diameter-bound}
\end{equation}
this implies $V(C)\leq h^{d+1}$.
On the other hand, the mesh must be \defm{quasiuniform} in the following sense:
\begin{equation}
    V(C)=\Omega(h^{d+1})\qquad(C\in\sys C^h). \label{eq:quasi-uniformity}
\end{equation}
Moreover, the cell surface measure must be controlled:
\begin{equation}
    S(\partial C)=O(h^d)\qquad(C\in\sys C^h). \label{eq:cell-boundary-bounded-measure}
\end{equation}

%
%

Let $\sys B^h$ be the $\sigma$-algebra generated by $\sys C^h$ over $\RC$; the elements of $M(\sys B^h;P)$ (the space
of $\sys B^h$-Borel-measurable maps on $\RC$ into $P$)
are called \defm{grid functions}. 
For $u^h\in M(\sys B^h;P)$, let $u^h_C\in P$ denote the constant value of $u^h$ on $\interior C$ ($C\in\sys C^h$).


\subsection{Numerical fluxes} 
Over every face $F\in\hat{\sys F}^h$ there is a \defm{numerical flux} $\abb{\Eta_F}{M(\sys B^h)}{\R}$.
(Note: the following definitions make sense for numerical entropy fluxes. The usual numerical fluxes can be reduced to this case;
see Section \ref{section:example}.)
The following requirements are imposed on numerical fluxes over interior faces.
\begin{enumerate}
\item 
    \defm{Consistency}: For $w\in P$, let $\hat w$ be the constant grid function with value $w$ (i.e.\ $\hat w_C=w$ for all $C\in\sys C^h$).
    We require
    \begin{alignat}{1}
        \Eta_{C\rarrow N}(\hat w) = \int_{C\rarrow N}\Vec\eta(w,y)\cdot n_{C\rarrow N}(y)dS(y). \label{eq:consistency}
    \end{alignat}
\item 
    \defm{Uniform continuity}: 
    There is a function $\delta_\Eta:(0,\infty)\rightarrow(0,\infty)$ so that
    \begin{alignat}{1}
        & \forall~\epsilon>0,~h>0,~F\in\sys F^h,~w\in\R^m,~w^h\in L^\infty(\sys B^h;P) : \notag\\
        & \|w^h-\hat w\|_{L^\infty(\sys B^h;P)}\leq\delta_\Eta(\epsilon)\ \Rightarrow\ |\Eta_F(w^h)-\Eta_F(\hat w)|\leq\epsilon h^d \label{eq:uniform-continuity}
    \end{alignat}
    (again, $\hat w$ denotes the constant grid function with value $w$).
\item
    \defm{Uniform boundedness}: for any\footnote{To apply this to common cases, restrict $P$ to be bounded (scalar case), bounded away from vacuum (gas dynamics) etc.} sequence $(w^h)_{h>0}$, 
    \begin{alignat}{1}
        \Eta_F(w^h) = O(h^d)\qquad(F\in\sys F^h). \label{eq:uniform-boundedness}
    \end{alignat}
\item
    \defm{Conservativeness}: for all $w^h\in M(\sys B^h;P)$, $C\rarrow N\in\sys F^h$,
    \begin{equation}
        \Eta_{C\rarrow N}(w^h) = - \Eta_{N\rightarrow C}(w^h). \label{eq:conservation-property}
    \end{equation}
\item
    \defm{Bounded stencil}: define the \defm{stencil} of $F\in\sys F^h$ as
    \begin{equation}
        \stn F:=\{C\in\sys C^h:\text{$w^h\mapsto \Eta_F(w^h)$ not constant in $w^h_C$}\} \notag
    \end{equation}
    and require
    \begin{equation}
        \sup_{C\in\stn F}d(C,F)=O(h). \label{eq:bounded-stencil-1}
    \end{equation}
\end{enumerate}
For initial faces, we impose the \defm{numerical initial condition}
\begin{equation}
    \Eta_{\partial\rarrow C}(w^h) = -\Eta_{C\rarrow\partial}(w^h) = \int_{\partial\rarrow C}\eta_0(u_0(x),0,x)dS(x)\qquad\forall w^h\in M(\sys B^h;\R^m). 
    \label{eq:inicond}
\end{equation}

\subsection{Numerical sources}
\label{subsection:numerical-sources}

The source terms in \eqref{weak-entropy-condition} are approximated by \defm{numerical sources}: functions
$G_C:M(\sys B^h;P)\rightarrow\R$ for each cell $C\in\sys C^h$. The numerical sources must satisfy the following conditions (that
are very similar to the ones for numerical fluxes):
\begin{enumerate}
\item 
    \defm{Consistency}:
    \begin{alignat}{1}
        \forall w\in P,~C\in\sys C^h:G_C(\hat w) = \int_Cg(w,y)dV(y) \label{eq:consistency-sources}
    \end{alignat}
    (where $\hat w$ is the constant grid function with value $w\in\R^m$). 
\item 
    \defm{Uniform continuity}: there is a function $\delta_G:(0,\infty)\rightarrow(0,\infty)$ so that
    \begin{alignat}{1}
        & \forall~\epsilon>0,~h>0,~C\in\sys C^h,~w\in P,~w^h\in L^\infty(\sys B^h;P) : \notag \\
        & \|w^h-\hat w\|_{L^\infty(\sys B^h;P)}\leq\delta_G(\epsilon)\ \Rightarrow\ |G_C(w^h)-G_C(\hat w)|\leq\epsilon h^{d+1}. \label{eq:uniform-continuity-sources}
    \end{alignat}
\item
    \defm{Uniform boundedness}: for any sequence $(w^h)_{h>0}$, 
    \begin{alignat}{1}
        G_C(w^h) = O(h^{d+1})\qquad(C\in\sys C^h). \label{eq:uniform-boundedness-sources}
    \end{alignat}
\item
    \defm{Bounded stencil}: define the \defm{stencil} of $C\in\sys C^h$ as
    \begin{equation}
        \stn C:=\{C'\in\sys C^h:\text{$w^h\mapsto G_C(w^h)$ not constant in $w^h_{C'}$}\} \notag
    \end{equation}
    and require
    \begin{equation}
        \sup_{C'\in\stn C}d(C,C')=O(h). \label{eq:bounded-stencil-sources}
    \end{equation}
\end{enumerate}

\subsection{Result}

\begin{thm}
    \label{thm:lax-wendroff-ich}%
    If a sequence $(u^h)_{h>0}$ of grid functions satisfies the discrete scalar inequalities
    \begin{alignat}{1}
        \sum_{N\in\hat{\sys C}^h,~C\cap N\neq\emptyset}\Eta_{C\rarrow N}(u^h) \leq G_C(u^h) \qquad (\forall h>0,~C\in\sys C^h) \label{eq:numerical}
    \end{alignat}
    and converges almost everywhere to $u$, then $u$ satisfies \eqref{weak-entropy-condition}. 
\end{thm}

\subsection{An example}

\label{section:example}

For illustration, the Lax-Friedrichs scheme (for a system $u_t+f(u)_x=0$ with initial condition $u(0,x)=u_0(x)$ on a uniform 1D 
grid with cell size $h$ and uniform time steps $\lambda h$ ($0<\lambda\leq 1$ constant))
is fit into the abstract framework in the previous sections.
For $h>0$, let $\sys C^h$ contain the cells $C_j^n$,
\begin{alignat}{1}
    C_j^n:=[jh,(j+1)h]\times[n\lambda h,(n+1)\lambda h] \qquad(j\in\Z,\ n\in\Nnull).
\end{alignat}
Numerical fluxes: for $j\in\Z$, $n\in\Nnull$,
\begin{alignat}{1}
    F_{\partial\rarrow C_j^0}(w^h) &:= \int_0^hu_0(jh+y)~dy, \\
    F_{C_j^n\rarrow C_j^{n+1}}(w^h) &:= hw^h_{C_j^{n+1}}, \\
    F_{C_j^n\rarrow C_{j+1}^n}(w^h) &:= h\left(\lambda\frac{f(w^h_{C_j^n})+f(w^h_{C_{j+1}^n})}{2} - \frac{w^h_{C_{j+1}^n}-w^h_{C_j^n}}{2}\right)
\end{alignat}
Numerical sources: all $=0$.
It is easy to check the numerical fluxes satisfy all conditions, in particular consistency. The numerical solutions $u^h\in M(\sys B^d;\R^m)$ 
are defined by
\begin{alignat}{1}
    C_0^j &:= h^{-1}\int_{jh}^{(j+1)h}u_0(x)dS(x), \\
    \sum_{N\in\hat{\sys C}^h}F_{N\rarrow C}(u^h) &= 0 \qquad (\forall h>0,~\forall C\in\sys C^h) \label{laxf}
\end{alignat}
which is exactly the literature definition, using different notation.

\eqref{laxf} is a system of $\R^m$-valued equations for $u^h$,
but it can obviously be converted into $2\cdot m$ systems of scalar inequalities of the type \eqref{eq:numerical}; 
Theorem \ref{thm:lax-wendroff-ich} applied to each of them separately implies \eqref{weak-conservation-law}, i.e.\ that $u$ is a weak solution.

In a similar fashion, it can be verified that the limit is an entropy solution. For Burgers equation ($f(u)=\ef{2}u^2$), 
it is sufficient to prove the entropy inequality 
for the Kru\v{z}kov family of entropies and entropy fluxes, 
\begin{alignat}{1}
    \eta_0(u) = |u-a|, \qquad \eta_1(u) = \operatorname{sgn}(u-a)f(u),
\end{alignat}
where $a\in\R$ is the family parameter. The numerical entropy fluxes
\begin{alignat}{1}
    \Eta_{\partial\rarrow C_j^0}(u^h) &:= \int_0^h\eta_0(u_0(jh+y))~dy, \\
    \Eta_{C_j^n\rarrow C_j^{n+1}}(u^h) &:= h\eta_0(u^h_{C_j^n}), \\
    \Eta_{C_j^n\rarrow C_{j+1}^n}(u^h) &:= \lambda h(F_{C_j^n\rarrow C_{j+1}^n}(u^h\vee c) - F_{C_j^n\rarrow C_{j+1}^n}(u^h\wedge c))
\end{alignat}
(see \cite{crandall-majda}) are consistent with $\Vec\eta$ and satisfy the discrete entropy inequality \eqref{eq:numerical}, so Theorem \ref{thm:lax-wendroff-ich} 
asserts that $u$ is an entropy solution (in the sense \eqref{weak-entropy-condition}).

\section{Proof of Theorem \ref{thm:lax-wendroff-ich}}

\label{section:fully-discrete}

The proof is based on two essential ideas: 
first, the original proof in \cite{lax-wendroff} uses summation by parts (in analogy to the integration by parts
used to derive the concept of weak solution); this requires a Cartesian grid. This obstacle is bypassed by \emph{approximating}
cubes with sidelength $H$ in a uniform Cartesian grid by cells in an unstructured grid with parameter $h$; see Figure \ref{grids}. 
Summation by parts is carried out for these cubes.

Second, since $(u^h)$ converges in $L^1_{\loc}(\RC)$, $u$ and $u^h$ will be ``close'' and ``nearly constant'' in a suitable 
neighbourhood of ``almost all'' cubes (for $h\downarrow 0$), so the continuity and consistency properties of numerical fluxes and sources 
can be exploited. For the ``few'' remaining ``bad'' cubes, one can use uniform boundedness of numerical fluxes and sources.
The proof will be completed by first fixing a sufficiently small ratio $\rho$
to minimize geometric errors and then choosing a sufficiently small $H>0$ to control integral errors. 

\subsection{Cubes}
Let $e^{(i)}\in\Z^{d+1}$ (with $i\in\{0,\dotsc,d\}$) be
the standard basis vectors, with $i$th component $=1$, all other components $=0$.
Omitting the parameter $H$ for readability, define the closed cubes
\begin{alignat}{1}
    I_{k} &:= H\cdot\prod_{i=0}^d[k_i,k_i+1]\qquad(k\in\Nnull\times\Z^d)\notag
\end{alignat}
with faces
\begin{alignat}{1}
    \partial I_{k}^{i+} &:= H\cdot\left(\prod_{j=0}^{i-1}[k_j,k_j+1]\times\{k_i+1\}\times\prod_{j=i+1}^d[k_j,k_j+1]\right), \notag\\
    \partial I_{k}^{i-} &:= H\cdot\left(\prod_{j=0}^{i-1}[k_j,k_j+1]\times\{k_i\}\times\prod_{j=i+1}^d[k_j,k_j+1]\right); \notag
\end{alignat}
note that the interiors of the $I_{k}$ are pairwise disjoint and that 
\begin{alignat}{1}
    \overline{\RC} = \bigcup_{k\in\ZC}I_{k};
\end{alignat}
moreover
\begin{alignat}{1}
    \partial I_{k} = \bigcup_{i=0}^d(\partial I_{k}^{i-}\cup\partial I_{k}^{i+})\qquad(k\in\Nnull\times\Z^d).
\end{alignat}

\subsection{Cube approximation}
For given $H,h>0$ ($h<H/2$) and $k\in\ZC$, we select an approximation $\tilde I_{k}\subset\sys C^h$ to $I_{k}$ by requiring
that 
\begin{alignat}{1}
    C\cap I_{k}\neq\emptyset\qquad\forall C\in\tilde I_{k} \label{eq:tildeQk}
\end{alignat}
\emph{and} that the sets $\tilde I_{k}$ form
a partition of $\sys C^h$. (These conditions need not determine $\tilde I_{k}$ uniquely; the particular choice is not important.
\eqref{eq:tildeQk} admits the existence of such $\tilde I_{k}$ because the $I_{k}$ cover $\overline{\RC}$.)

\begin{figure}
\input{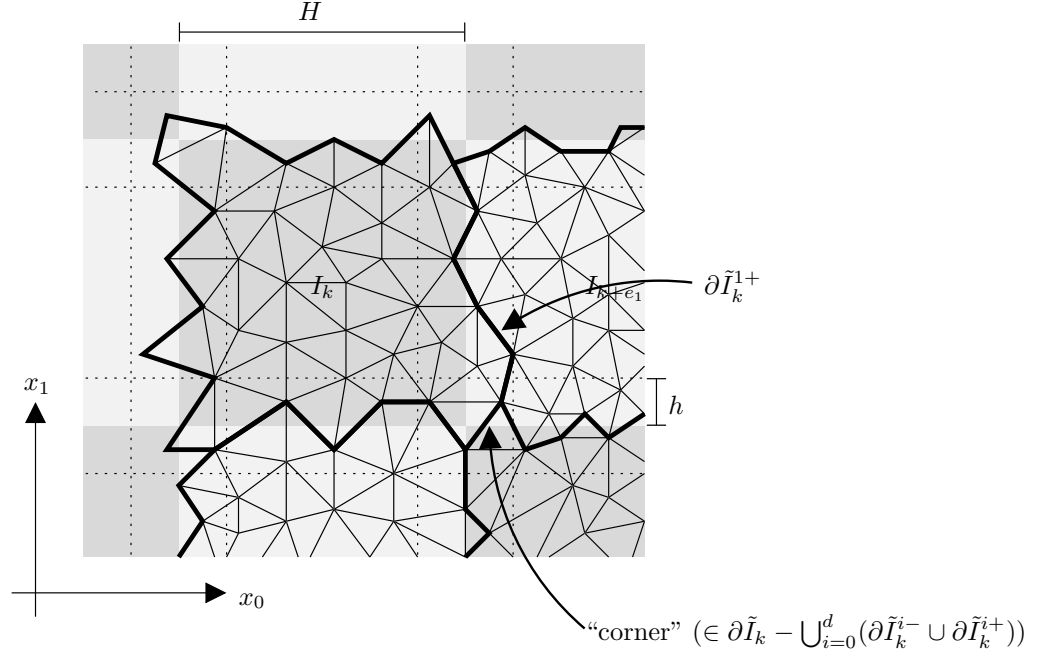}
\caption{A cube $I_{k}$ (checkerboard-shaded grid) is approximated by a cluster (thick boundary) of grid cells (thin triangles).}
\label{grids}
\end{figure}

We also define (for $i\in\Nd$, $k\in\Nnull\times\Z^d$)
\begin{alignat}{1}
    \partial\tilde I_{k} &:= \{C\rarrow\partial\in\hat{\sys F}^h:C\in\tilde I_{k}\}\cup\{C\rarrow N\in\sys F^h:C\in\tilde I_{k},~N\not\in\tilde I_{k}\}.
    \label{eq:partialtildeQk-2} \\
    \partial\tilde I_{k}^{i\pm} &:= 
    \begin{cases}
        \{C\rarrow N\in\sys F^h:C\in\tilde I_{k},\ N\in\tilde I_{k\pm e^{(i)}}\}, & k\pm e^{(i)}\in\ZC \\
        \{C\rarrow\partial\in\hat{\sys F}^h:C\in\tilde I_{(0,k_1,\dotsc,k_d)}\}, & \text{else} \\
    \end{cases},  \label{eq:partialtildeQk-1} 
\end{alignat}
Note that 
\[ 
C\rarrow N\in \partial\tilde I_{k+e^{(i)}}^{i-}\qquad\Leftrightarrow\qquad N\rarrow C\in\partial\tilde I_{k}^{i+};
\]
however
\[ \bigcup_{i=0}^d(\partial\tilde I_{k}^{i-}\cup\partial\tilde I_{k}^{i+})= \partial\tilde I_{k} \]
is \emph{not} true in general (see Lemma \ref{lemma:diagonal-faces}) because some faces belong to ``corners'' rather than sides of the approximated 
cubes (see Figure \ref{grids}).


\subsection{Geometric estimates}

For $r>0$ and $A\subset\overline{\RC}$, define the closed neighbourhoods
\begin{alignat}{1}
    \overline Q_r(A)&:=\{y\in\overline{\RC}:d(y,A)\leq r\} \notag 
\end{alignat}

\begin{lemma}
    \begin{alignat}{1}
        \bigcup_{C\in\tilde I_{k}}C &\subset \overline Q_h(I_{k}),\qquad I_{k} \subset \overline Q_h(\bigcup_{C\in\tilde I_{k}}C), \label{eq:Ik-tildeIk}
    \end{alignat}
    so
    \begin{alignat}{1}
        V(\bigcup_{C\in\tilde I_{k}}C) &= H^{d+1}+O(\rho H^{d+1}). \label{eq:VtildeIk}
    \end{alignat}
    Moreover
    \begin{alignat}{1}
        \bigcup_{F\in\partial\tilde I^{i\pm}_{k}}F &\subset \overline Q_h(\partial I^{i\pm}_{k}) \label{eq:partialtildeIkis-partialIkis}.
    \end{alignat}
\end{lemma}
\begin{proof}
    These are immediate consequences of \eqref{eq:diameter-bound} and of \eqref{eq:tildeQk}, \eqref{eq:partialtildeQk-2} resp.\ \eqref{eq:partialtildeQk-1}. 
\end{proof}


\begin{lemma}
    \label{large-small-bounded}%
    $A\subset\RC$ meets $\leq\frac{V\left(\overline Q_h(A)\right)}{h^{d+1}}$ cells in $\sys C^h$.
\end{lemma}
\begin{proof}
    $\diam C\leq h$, so $C\cap A\neq\emptyset$ implies $C\subset\overline Q_h(A)$. 
    However, \eqref{eq:quasi-uniformity} implies that
    $\overline Q_h(A)$ cannot contain more than $V(\overline Q_h(A))h^{-(d+1)}$ cells.
\end{proof}

\begin{cor}
    \begin{alignat}{1}
        \sup_{F\in\hat{\sys F}^h}\#\stn F&=O(1); \label{eq:cardinality-stencil} \\
        \sup_{C\in\sys C^h}\#\{N\in\sys C^h:C\cap N\neq\emptyset\} &= O(1); \label{eq:cardinality-faces}
    \end{alignat}
    for all $k\in\ZC$ and $i\in\{0,\dotsc,d\}$,
    \begin{alignat}{1}
        \sup_{k\in\Nnull\times\Z^d}\#\tilde I_{k} &= O(\rho^{-d-1}), \label{eq:cardinality-tildeQk} \\
        \sup_{k\in\Nnull\times\Z^d,i\in\{0,\dotsc,d\},s\in\{+,-\}}\#\partial\tilde I_{k}^{is} &= O(\rho^{-d}). \label{eq:cardinality-partialtildeQk}
    \end{alignat}
\end{cor}
\begin{proof}
    \eqref{eq:cardinality-stencil}: by \eqref{eq:bounded-stencil-1}, there is a constant $c$ (independent of $h$) so that
    \[ \bigcup_{C\in\stn F}C\subset\overline Q_{ch}(F) \]
    for all faces $F$. Since $\diam F\leq h$, $\overline Q_{ch}(F)$ is contained in some closed ball with diameter $O(h)$, hence volume $O(h^{d+1})$.
    At most $O(1)$ cells fit into this ball, so $\#\stn F=O(1)$.

    \eqref{eq:cardinality-faces}: this is immediate from Lemma \ref{large-small-bounded}.

    \eqref{eq:cardinality-tildeQk}: \eqref{eq:Ik-tildeIk} implies
    \[
    \overline Q_h(\bigcup_{C\in\tilde I_{k}}C) \subset \overline Q_{2h}(I_{k}),
    \]
    hence ($\rho\leq\frac{1}{2}$)
    \[
    V(\overline Q_h(\bigcup_{C\in\tilde I_{k}}C)) = O(H^{d+1}).
    \]
    Hence by Lemma \ref{large-small-bounded}, at most $O(\frac{H^{d+1}}{h^{d+1}})=O(\rho^{-d-1})$ cells meet $\bigcup_{C\in\tilde I_{k}}C$.

    \eqref{eq:cardinality-partialtildeQk}: from \eqref{eq:partialtildeIkis-partialIkis} derive
    \[
    \overline Q_h(\bigcup_{F\in\partial\tilde I^{i\pm}_{k}}F)\subset\overline Q_{2h}(\partial I^{i\pm}_{k}),
    \]
    so
    \[
    V(\overline Q_h(\bigcup_{F\in\partial\tilde I^{i\pm}_{k}}F)) = O(hH^d) = O(\rho H^{d+1}).
    \]
    Hence Lemma \ref{large-small-bounded} shows that at most $O(\frac{\rho H^{d+1}}{h^{d+1}})=O(\rho^{-d})$ cells 
    meet $\bigcup_{F\in\partial\tilde I^{i\pm}_{k}}F$; by \eqref{eq:cardinality-faces} each has $O(1)$ faces.
\end{proof}

The following lemma states that, for small $\rho$, ``most'' of $\partial\tilde I_{k}$ is composed of the $\partial\tilde I_{k}^{i\pm}$ ($i=0,\dotsc,d$),
i.e.\ we can ignore the ``corners'' of $I_{k}$.

\begin{figure}
\input{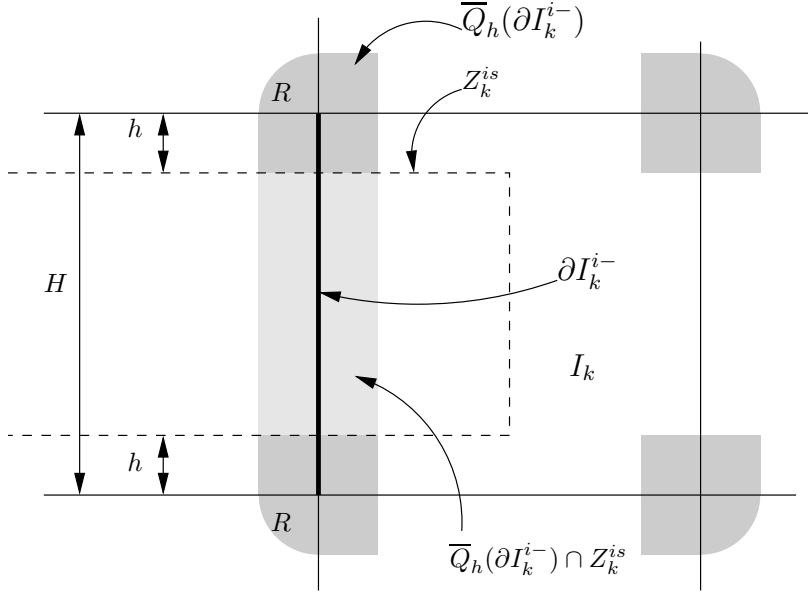}
\caption{$V(R)=O(h^2H^{d-1})=O(\rho^2H^{d+1})$}
\label{R}
\end{figure}

\begin{lemma}
    \label{lemma:diagonal-faces}%
    Define the ``half-cylinders''
    \begin{alignat}{1}
        Z^{i\pm}_{k} &:= \{y\in\overline{\RC}:y_i\gtrless H(k_i+\ef{2}),\ y_j\in H\cdot(k_j+\rho,k_j+1-\rho)\quad (j\neq i)\} 
    \end{alignat}
    (see Figure \ref{R}). Then
    \begin{alignat}{1}
        \sum_{F\in\partial\tilde I^{i\pm}_{k}}S(F) = \sum_{F\in\partial\tilde I^{i\pm}_{k}}S(F\cap Z^{i\pm}_{k}) + O(\rho)H^d. \label{eq:corners-1}
    \end{alignat}
    Moreover,
    \begin{alignat}{1}
        \sum_{F\in\partial\tilde I_{k}-\bigcup_{i=0}^d(\partial\tilde I^{i-}_{k}\cup\partial\tilde I^{i+}_{k})}S(F) = O(\rho)H^d.
        \label{eq:corners-2}
    \end{alignat}
\end{lemma}
\begin{proof}
    (See Figure \ref{R}.)
    Let $F=C\rarrow N\in\partial\tilde I_{k}^{is}$ for some $i\in\{0,\dotsc,d\},s\in\{+,-\}$ (or $F=C\rarrow\partial$).
    By \eqref{eq:partialtildeIkis-partialIkis}, 
    whenever $F\cap Z^{is}_{k}\neq F$, then $F$ (and hence $C$) meets 
    $$R:=\overline Q_h(\partial I_{k})-\bigcup_{i=0}^d(Z^{i+}_{k}\cup Z^{i-}_{k}).$$
    However, it is easy to verify that
    \begin{alignat}{1}
        V(\overline Q_h(R)) = O(h^2H^{d-1}) = O(\rho^2H^{d+1}).
    \end{alignat}
    By Lemma \ref{large-small-bounded}, at most
    \begin{alignat}{1}
        O(\frac{V(\overline Q_h(R))}{h^{d+1}})&=O(\rho^{1-d}) \notag
    \end{alignat}
    cells can meet $R$. By \eqref{eq:cell-boundary-bounded-measure}, their total surface measure is 
    \begin{alignat}{1}
        =O(\rho^{1-d})O(h^d)=O(\rho H^d); \notag
    \end{alignat}
    this implies \eqref{eq:corners-1}.
    
    Regarding \eqref{eq:corners-2}: 
    whenever $C\rarrow N\in\partial\tilde I_{k}$ (with $C\in\tilde I_{k}$) is not contained in any $\partial\tilde I_{k}^{i\pm}$, 
    then $C\rarrow N$ belongs to a ``corner'', i.e.\ $N\in\tilde I_{k+m}$ for some $m\in\Z^{d+1}$ with $|m|_\infty=1$, $|m|_1\geq 2$ (because of
    $\rho<\ef{2}$ and $\diam C,\diam N\leq h$).
    This means $C\rarrow N\subset R$; \eqref{eq:corners-1} states
    that these faces may be ignored at a cost of surface measure of $O(\rho H^d)$. 

    $C\rarrow\partial\in\partial\tilde I^{0-}_{k}$ for some $k\in\{0\}\times\Z^d$, so initial faces do not contribute to the sum in
    \eqref{eq:corners-2}.
\end{proof}

The numerical fluxes over $\bigcup_{F\in\partial\tilde I_{k}^{i\pm}}F$ will be pieced together to approximate the exact flux
over $\partial I_{k}^{i\pm}$. This requires the following geometric estimate.

\begin{figure}
\input{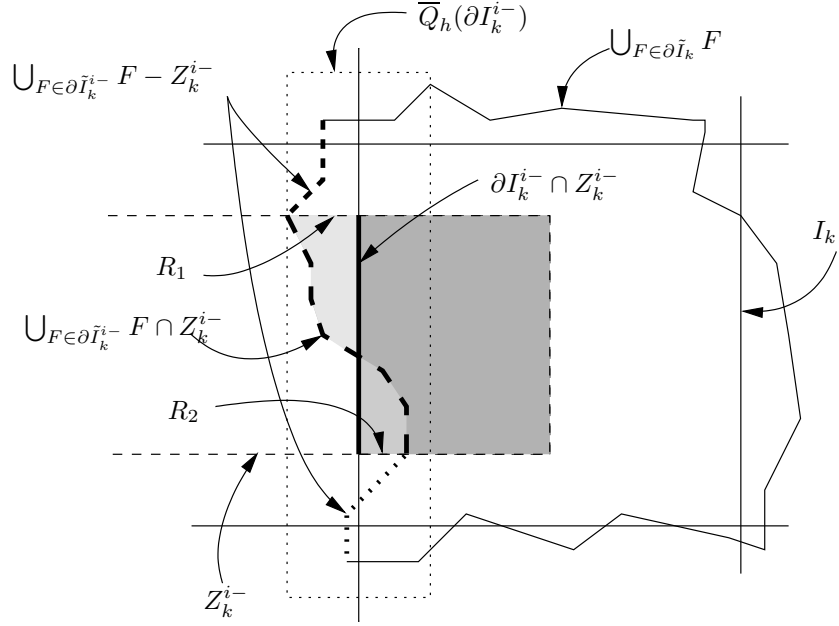}
\caption{The surfaces $R_1$, $R_2$, $\bigcup_{F\in\partial\tilde I_{k}^{i-}}F-Z_{k}^{i-}$, and $\partial I_{k}^{i-}-Z_{k}^{i-}$ have measure $O(\rho)H^d$
    and can be neglected.} 
\label{fig:nint}
\end{figure}

\begin{lemma}
    \label{lemma:nint}%
    For all $k\in\ZC$, $i\in\Nd$,
    \begin{alignat}{1}
        \left|\left(\int_{\bigcup_{F\in\partial\tilde I^{i-}_{k}}F}-\int_{\partial I^{i-}_{k}}\right)\Vec n~dS\right| &=O(\rho)H^d. 
        \label{eq:nint}
    \end{alignat}
\end{lemma}
\begin{proof}
    (See Figure \ref{fig:nint}.)
    \begin{alignat}{1}
        & \left(\int_{\bigcup_{F\in\partial\tilde I^{i-}_{k}}F}-\int_{\partial I^{i-}_{k}}\right)n~dS \notag \\
        &\topanno{\eqref{eq:corners-2}}{=} 
        \left(\int_{\bigcup_{F\in\partial\tilde I^{i-}_{k}}F\cap Z^{i-}_{k}}-\int_{\partial I^{i-}_{k}\cap Z^{i-}_{k}}\right)n~dS + O(\rho)H^d \label{eq:nint-1}
    \end{alignat}
    The first summand in \eqref{eq:nint-1} equals
    \begin{alignat}{1}
        &\left(\int_{\partial(\bigcup_{C\in\tilde I_{k}C}\cap Z^{i-}_{k})}-\int_{\partial(I_{k}\cap Z^{i-}_{k})}\right)n~dS  \notag \\
        &=0 \notag
    \end{alignat}
    up to 
    \begin{alignat}{1}
        \int_{R_1}n dS - \int_{R_2}n dS \label{eq:rest}
    \end{alignat}
    where $R_1,R_2$ are contained in
    \begin{alignat}{1}
        \partial Z^{i-}_{k}\cap\overline Q_h(\partial I^{i-}_{k}),\notag
    \end{alignat}
    a set with total surface measure $O(\rho)H^d$; hence the two integrals in \eqref{eq:rest} are $O(\rho)H^d$ too.
\end{proof}

\subsection{Completion of the proof}

Consider an arbitrary nonnegative test function $\phi\in C_c^\infty(\overline\RC)$.
Since $\phi$ has compact support, it is sufficient to consider the finite subsets
\begin{alignat}{1}
    K &:= \{k\in\ZC:\exists\ell\in\Z^{d+1}:|\ell|\leq 1,~\supp\phi\cap I_{k+\ell}\neq\emptyset\}\notag
\end{alignat}
of cube indices. Note that
$$\#K=O(\rho^{-(d+1)}).$$

Define
\begin{alignat}{1}
    \stn k &:= \bigcup_{F\in\partial\tilde I_{k}}\stn F\cup\bigcup_{C\in\tilde I_{k}}\stn C, \label{eq:stn-k} \\
    \stn K &:= \bigcup_{k\in K}\stn k \label{eq:stn-K}.
\end{alignat}

\begin{lemma}
    For all $w\in L^1_{\loc}(\overline{\RC})$,
    \begin{alignat}{1}
        \sum_{k\in K}\int_{\bigcup_{C\in\stn k}C}|w(y)|dV(y) & = O(1)\int_{\bigcup_{C\in\stn K}C}|w(y)|dV(y) \label{eq:stn-Qk-comparison}
    \end{alignat}
\end{lemma}
\begin{proof}
    Due to bounded stencils \eqref{eq:bounded-stencil-1} resp.\ \eqref{eq:bounded-stencil-sources}, 
    bounded diameters \eqref{eq:diameter-bound} and $\rho\leq\ef{2}$ (see \eqref{eq:rho-def}), 
    \begin{alignat*}{1}
        \sup_{C\in\sys C^h}\#\{k\in K:C\in\stn k\} &= O(1).
    \end{alignat*}
    Hence
    \begin{alignat*}{1}
        \sum_{k\in K}\int_{\bigcup_{C\in\stn k}C}|w(y)|dV(y) &= \sum_{C\in\stn K}\#\{k\in K:C\in\stn k\}\int_C|w(y)|dV(y) \\
        &\leq O(1)\sum_{C\in\stn K}\int_C|w(y)|dV(y).
    \end{alignat*}
\end{proof}

We need to show that $u,u^h$ are ``almost constant'' and ``close'' on ``most'' cubes. We introduce a new parameter $\epsilon>0$.
Again omitting $\epsilon,\rho,H>0$ from the symbols for readability, define 
\begin{alignat}{1}
    \overline u_{k} &:= V(I_{k})^{-1}\int_{I_{k}}u(y)dV(y), \\
    B_1 &:= \{k\in K:\exists~C\in\stn k:\ |u^h_C-\overline u_{k}|\text{$>\min\{\delta_\Eta(\epsilon),\delta_G(\epsilon)\}$}\}, \label{eq:B1def} \\
    B_2 &:= \{k\in K:\frac{V\{y\in I_{k}:|u(y)-\overline u_{k}|\text{$>\min\{\delta_\Eta(\epsilon),\delta_G(\epsilon)\}$}\}}{V(I_{k})}>\epsilon \}, \label{eq:B2def} \\
    B &:= B_1\cup B_2, \notag \\
    G &:= K-B \label{eq:Gdef};
\end{alignat}
$B$ contains the ``bad'', $G$ the ``good'' cube indices.

\begin{lemma}
    \label{th:Qk-nbh}%
    For any choice of $\epsilon,\rho>0$
    \begin{alignat}{1}
        \frac{\#B}{H^{-(d+1)}} &= o_{\rho,\epsilon}(1). \label{eq:Qk-nbh}
    \end{alignat}
\end{lemma}
\begin{proof}
    $B_1$ and $B_2$ are treated separately. First $B_2$:
    let $w\in C^\infty(\overline{\RC})$ be arbitrary (it will approximate $u$); define 
    \begin{alignat}{1}
        \overline w_{k} &= V(I_{k})^{-1}\int_{I_{k}}w(y)dV(y). \notag
    \end{alignat}
    Then
    \begin{alignat}{1}
        & \sum_{k\in B_2}\int_{I_{k}}|u(y)-\overline u_{k}|dV(y) \notag \\
        &\leq \sum_{k\in K}\int_{I_{k}}|u(y)-\overline u_{k}|dV(y) \notag \\
        &\leq \sum_{k\in K}\int_{I_{k}}|u(y)-w(y)|dV(y) 
        + \sum_{k\in K}\int_{I_{k}}|w(y)-\overline w_{k}|dV(y) 
        + \sum_{k\in K}V(I_{k})|\overline w_{k}-\overline u_{k}| \notag \\
        & = O(\|u-w\|_{L^1(\bigcup_{k\in K}I_{k})} + \|Dw\|_{L^\infty(\bigcup_{k\in K}I_{k})}H) \label{eq:B2-1} 
    \end{alignat}
    because
    \begin{alignat}{1}
        & \sum_{k\in K}V(I_{k})|\overline u_{k}-\overline w_{k}| = \sum_{k\in K}\left|\int_{I_{k}}u(y)-w(y)dV(y)\right| \notag \\
        & \leq \sum_{k\in K}\int_{I_{k}}|u(y)-w(y)|dV(y) \leq \|u-w\|_{L^1(\bigcup_{k\in K}I_{k})} \notag
    \end{alignat}
    and because
    \begin{alignat}{1}
        |\overline w_{k}-w(y)| &= O(H\|Dw\|_{L^\infty(\bigcup_{k\in K}I_{k})}) \notag
    \end{alignat}
    for $y\in I_{k}$, by smoothness of $w$.
    For each $k\in B_2$,
    \begin{alignat}{1}
        \int_{I_{k}}|u(y)-\overline u_{k}|dV(y) &\geq \epsilon\min\{\delta_\Eta(\epsilon),\delta_G(\epsilon)\}V(I_{k}),
    \end{alignat}
    by definition \eqref{eq:B2def} of $B_2$, so \eqref{eq:B2-1} implies that
    \begin{alignat}{1}
        \#B_2\cdot H^{d+1} = \sum_{k\in B_2}V(I_{k}) &= O(\frac{\|u-w\|_{L^1(\bigcup_{k\in K}I_{k})}+\|Dw\|_{L^\infty(\bigcup_{k\in K}I_{k})}H}{\epsilon\min\{\delta_\Eta(\epsilon),\delta_G(\epsilon)\}}).
    \end{alignat}
    By first choosing $w\in C^\infty(\RC)$ with sufficiently small $\|u-w\|_{L^1(\bigcup_{k\in K}I_{k})}$ 
    and then choosing an upper bound for $H$, the right-hand side can be made arbitrarily small.

    Regarding $B_1$,
    \begin{alignat}{1}
        & \sum_{k\in B_1}\int_{\bigcup_{C\in\stn k}C}|u^h(y)-\overline u_{k}|dV(y) \notag \\
        & \leq \sum_{k\in K}\int_{\bigcup_{C\in\stn k}C}|u^h(y)-\overline u_{k}|dV(y) \notag \\
        & \leq \sum_{k\in K}\int_{\bigcup_{C\in\stn k}C}|u^h(y)-u(y)|dV(y)
        + \sum_{k\in K}\int_{\bigcup_{C\in\stn k}C}|u(y)-w(y)|dV(y) \notag \\
        & + \sum_{k\in K}\int_{\bigcup_{C\in\stn k}C}|w(y)-\overline w_{k}|dV(y)
        + \sum_{k\in K}\int_{\bigcup_{C\in\stn k}C}|\overline w_{k}-\overline u_{k}|dV(y) \notag \\
        & \topanno{\eqref{eq:stn-Qk-comparison}}{=} O(\|u^h-u\|_{L^1(\bigcup_{C\in\stn K}C)}+\|u-w\|_{L^1(\bigcup_{C\in\stn K}C)} + H\|Dw\|_{L^\infty(\bigcup_{C\in\stn K}C)}). \label{eq:B1-1}
    \end{alignat}
    For each $k\in B_1$,
    \begin{alignat}{1}
        \int_{\bigcup_{C\in\stn k}C}|u^h(y)-\overline u_{k}|dV(y) &= \Omega(\rho^{d+1}H^{d+1}\min\{\delta_\Eta(\epsilon),\delta_G(\epsilon)\}) = \Omega(\rho^{d+1}V(I_{k})\min\{\delta_\Eta(\epsilon),\delta_G(\epsilon)\}), \notag
    \end{alignat}
    so \eqref{eq:B1-1} (with \eqref{eq:quasi-uniformity}) shows
    \begin{alignat}{1}
        \#B_1\cdot H^{d+1} &= O\left(\frac{\|u^h-u\|_{L^1(\bigcup_{C\in\stn K}C)}+\|u-w\|_{L^1(\bigcup_{C\in\stn K}C)}
            + \|Dw\|_{L^\infty(\bigcup_{C\in\stn K}C)}H}{\rho^{d+1}\min\{\delta_\Eta(\epsilon),\delta_G(\epsilon)\}}\right) \notag
    \end{alignat}
    By first choosing a suitable $w$ and then choosing an upper bound for $H$, the right-hand side can be made arbitrarily small.
\end{proof}

Now we can finish the proof of Theorem \ref{thm:lax-wendroff-ich}. 
Sum \eqref{eq:numerical} over $C\in\tilde I_{k}$, multiply with $\phi(Hk)$ ($\geq 0$) and sum over $k\in K$:
\begin{alignat}{1}
    \sum_{k\in K}\phi(Hk)\sum_{C\in\tilde I_{k}}G_C(u^h)
    &\geq
    \sum_{k\in K}\phi(Hk)\sum_{F\in\partial\tilde I_{k}}\Eta_F(u^h) \label{eq:bigsum7}
\end{alignat}
(here we used the conservation property \eqref{eq:conservation-property} to eliminate $F\not\in\bigcup_{k\in K}\partial\tilde I_{k}$,
i.e.\ $F=C\rarrow N$ with $C,N\in\tilde I_{k}$ for the same $k$).


Collecting the terms on the right-hand side of \eqref{eq:bigsum7} where $F$ is an initial face yields
\begin{alignat*}{1}
    \sum_{k\in K}\sum_{C\rarrow\partial\in\partial\tilde I^{0-}_{k}}\Eta_{C\rarrow\partial}(u^h)\phi(Hk) 
\end{alignat*}
which equals
\begin{alignat}{1}
    -\int_{\Rd}\eta_0(u_0(x))\phi(0,x)dS(x) + O(H), \label{eq:initial-faces}
\end{alignat}
using $\eqref{eq:inicond}$ and smoothness plus compact support of $\phi$. 

There are at most two terms per interior face in \eqref{eq:bigsum7}; they can be written (using \eqref{eq:conservation-property}) as
\begin{alignat}{1}
    \Eta_{C\rarrow N}(u^h)(\phi(Hk_C)-\phi(Hk_N)) \label{face-sums}
\end{alignat}
where $C\in \tilde I_{k_C}$, $N\in\tilde I_{k_N}$ ($|k_C-k_N|_\infty=1$).

Note that the $\phi$ difference is $O(H)$ (by smoothness of $\phi$), and
that for each $k\in K$, \eqref{eq:corners-2} allows to drop all terms for interior faces that do not belong to some $\partial\tilde I_{k}^{i\pm}$, 
at the cost of terms of size $O(\rho)H^d\cdot O(H)=O(\rho)H^{d+1}$ per $k$, hence $O(\rho)$ overall. Here, it is important that $\Eta_F(u^h)=O(h^d)$ 
(by \eqref{eq:uniform-boundedness}).

Moreover, by Lemma \ref{th:Qk-nbh} the number of bad cubes is $o_{\rho,\epsilon}(1)H^{-(d+1)}$, and due to uniform 
boundedness \eqref{eq:uniform-boundedness},
\eqref{eq:cell-boundary-bounded-measure} and \eqref{eq:cardinality-partialtildeQk},
\[ \sum_{C\rarrow N\in\partial\tilde I_{k}}\left|\Eta_{C\rarrow N}(u^h)\right| = O(H^d). \]
Since the $\phi$ difference supplies an extra $O(H)$, the sum of terms in \eqref{face-sums} where $k_C$ or $k_N$ is in $B$
is $o_{\rho,\epsilon}(1)H^{-(d+1)}\cdot O(H^d)\cdot O(H)=o_{\rho,\epsilon}(1)$. Hence the interior face part of \eqref{eq:bigsum7} is
\begin{alignat}{1}
    &=\sum_{i=0}^d\sum_{k\in G}\sum_{F\in\partial\tilde I^{i-}_{k+e^{(i)}}}\Eta_F(u^h)\subeq{(\phi(H(k+e^{(i)}))-\phi(Hk))}{}
    + O(\rho) + o_{\rho,\epsilon}(1) \notag \\
    &=
    \sum_{i=0}^d 
    \sum_{k\in G}
    \subeq{\sum_{F\in\partial\tilde I^{i-}_{k+e^{(i)}}}\Eta_F(u^h)}{}
    H^{-d}\int_{I_{k}}\frac{\partial\phi}{\partial y_i}(y)dV(y)
    + O(\rho)  + o_{\rho,\epsilon}(1)\notag 
    \intertext{(using smoothness and compact support of $\phi$)} \notag 
    &= \sum_{i=0}^d\sum_{k\in G}
    \subeq{\sum_{F\in\partial\tilde I^{i-}_{k+e^{(i)}}}\Eta_F(\hat{\overline u}_{k})}{}
    H^{-d}\int_{I_{k}}\frac{\partial\phi}{\partial y_i}(y)dV(y)
    + O(\rho+\epsilon)  + o_{\rho,\epsilon}(1)\notag 
    \intertext{(using the definition \eqref{eq:Gdef} of $G$ in uniform continuity \eqref{eq:uniform-continuity},  \notag
        combined with \eqref{eq:cardinality-partialtildeQk})} 
    &= \sum_{i=0}^d\sum_{k\in G}
    \subeq{\int_{\bigcup_{F\in\partial\tilde I^{i-}_{k+e^{(i)}}}F}\eta(\overline u_{k},\cdot)\cdot n~dS}{}
    H^{-d}\int_{I_{k}}\frac{\partial\phi}{\partial y_i}(y)dV(y)
    + O(\rho+\epsilon)  + o_{\rho,\epsilon}(1)\notag 
    \intertext{(using consistency \eqref{eq:consistency})} \notag 
    &= \sum_{i=0}^d\sum_{k\in G}
    \eta(\overline u_{k},Hk)\cdot\subeq{\int_{\bigcup_{F\in\partial\tilde I^{i-}_{k+e^{(i)}}}F}n~dS}{}~
    H^{-d}\int_{I_{k}}\frac{\partial\phi}{\partial y_i}(y)dV(y)
    + O(\rho+\epsilon)  + o_{\rho,\epsilon}(1)\notag 
    \intertext{(using smoothness in $y$ of $\eta$ (note \eqref{eq:partialtildeIkis-partialIkis},\eqref{eq:rho-def}) with \eqref{eq:cardinality-partialtildeQk} and \eqref{eq:cell-boundary-bounded-measure})} \notag 
    &\topanno{\eqref{eq:nint}}{=} -\sum_{i=0}^d\sum_{k\in G}\subeq{\eta_i(\overline u_{k},Hk) \int_{I_{k}}\frac{\partial\phi}{\partial y_i}dV(y)}{}
    + O(\epsilon+\rho) + o_{\rho,\epsilon}(1) \notag \\
    &= -\sum_{i=0}^d\subeq{\sum_{k\in G}}{}\int_{I_{k}}\frac{\partial\phi}{\partial y_i}(y)\eta_i(u(y),y)dV(y)  
    + O(\epsilon+\rho) + o_{\rho,\epsilon}(1) \notag 
    \intertext{(using the definition \eqref{eq:Gdef} of $G$ and smoothness in $y$ and $u$ of $\eta$)} \notag 
    &\topanno{\eqref{eq:Qk-nbh}}{=} -\sum_{i=0}^d\sum_{k\in K}\int_{I_{k}}\frac{\partial\phi}{\partial y_i}(y)\eta_i(u(y),y)dV(y)  
    + O(\epsilon+\rho) + o_{\rho,\epsilon}(1) \notag \\
    &= -\int_{\RC}\sum_{i=0}^d\frac{\partial\phi}{\partial y_i}(y)\eta_i(u(y),y)dV(y) 
    + O(\epsilon+\rho) + o_{\rho,\epsilon}(1). \label{eq:interior-faces}
\end{alignat}



It remains to treat the left-hand side of \eqref{eq:bigsum7}. As for the flux integrals, we may omit the terms for $k\in B$, 
at a cost of $o_{\rho,\epsilon}(1)H^{-(d+1)}\cdot O(H^{d+1})=o_{\rho,\epsilon}(1)$
(from \eqref{eq:Qk-nbh} resp.\ \eqref{eq:uniform-boundedness-sources}), hence:
\begin{alignat}{1}
    & \subeq{\sum_{k\in K}}{}\phi(Hk)\sum_{C\in\tilde I_{k}}G_C(u^h) \notag \\
    =& \sum_{k\in G}\phi(Hk)\subeq{\sum_{C\in\tilde I_{k}}G_C(u^h)}{}
    + o_{\rho,\epsilon}(1) \notag \\
    \topanno{\eqref{eq:uniform-continuity-sources}}{=} & \sum_{k\in G}\phi(Hk)\subeq{\sum_{C\in\tilde I_{k}}G_C(\hat{\overline u}_{k})}{} + O(\epsilon)
    + o_{\rho,\epsilon}(1) \notag \\
    \topanno{\eqref{eq:consistency-sources}}{=} & \sum_{k\in G}\phi(Hk)\subeq{\int_{\bigcup_{C\in\tilde I_{k}}C}g(\overline u_{k},y)dV(y)}{} + O(\epsilon)
    + o_{\rho,\epsilon}(1) \notag 
\end{alignat}
\clearpage
\begin{alignat}{1}
    =& \sum_{k\in G}\phi(Hk)\subeq{V\left(\bigcup_{C\in\tilde I_{k}}C\right)}{}g(\overline u_{k},Hk) + O(\epsilon)
    + o_{\rho,\epsilon}(1) \notag \\
    \topanno{\eqref{eq:VtildeIk}}{=} & \sum_{k\in G}\subeq{\phi(Hk)H^{d+1}}{}g(\overline u_{k},Hk)
    + O(\epsilon+\rho) + o_{\rho,\epsilon}(1) \notag \\
    =& \sum_{k\in G}\subeq{\int_{I_{k}}\phi(y)dV(y)g(\overline u_{k},Hk)}{}
    + O(\epsilon+\rho) + o_{\rho,\epsilon}(1) \notag \\
    =& \subeq{\sum_{k\in G}}{}\int_{I_{k}}\phi(y)g(u(y),y)dV(y) 
    + O(\epsilon+\rho) + o_{\rho,\epsilon}(1) \notag \\
    \topanno{\eqref{eq:Qk-nbh}}{=}& \subeq{\sum_{k\in K}\int_{I_{k}}}{}\phi(y)g(u(y),y)dV(y)
    + O(\epsilon+\rho) + o_{\rho,\epsilon}(1) \notag \\
    =& \int_{\RC}\phi(y)g(u(y),y)dV(y)
    + O(\epsilon+\rho) + o_{\rho,\epsilon}(1) \label{eq:bigsourcesum}
\end{alignat}

\subsection{Conclusion}

Combining \eqref{eq:initial-faces}, \eqref{eq:interior-faces} and \eqref{eq:bigsourcesum}, we get
\begin{alignat}{1}
    & \int_{\RC}\phi(y)g(u(y),y)dV(y) + O(\epsilon+\rho) + o_{\rho,\epsilon}(1) \notag \\
    \geq&
    -\int_{\RC}\sum_{i=0}^d\frac{\partial\phi}{\partial y_i}(y)\eta_i(u(y))dV(y) 
    - \int\eta_0(u_0(x))\phi(0,x)dS(x) \notag
\end{alignat}

Now, we can first make the $O(\epsilon+\rho)$ term arbitrarily small by choosing appropriate $\epsilon$ and $\rho$;
after that, the $o$ term can be made arbitrarily small as well 
by picking $H$. Therefore, $u$ satisfies \eqref{weak-entropy-condition}; the proof is complete.

\section{A counterexample for non-quasiuniform grids}

\label{section:counterexample}

\begin{figure}
\input{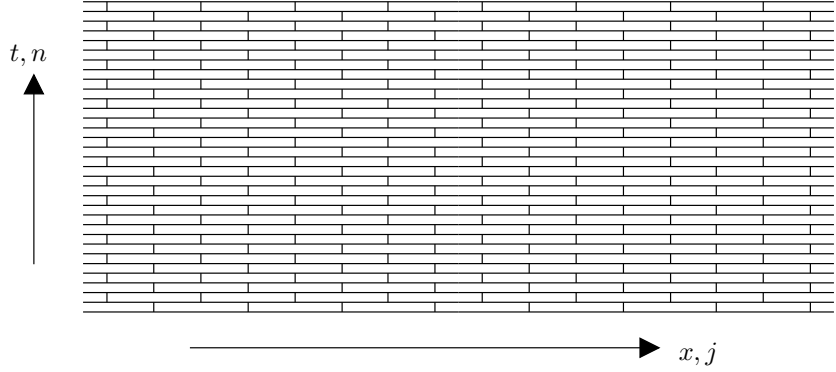}
\caption{Staggered Lax-Friedrichs with $\Delta t=h^3$, $\Delta x=h$: excessive refinement in $t$ direction causes oversmoothing and convergence
    to a non-solution.}
\label{fig:counterexample}
\end{figure}

While most assumptions in this paper are rather weak, an important exception is quasiuniformity 
(in the sense of \eqref{eq:quasi-uniformity}). Unfortunately, there is a strong counterexample to Theorem \ref{thm:lax-wendroff-ich} for 
non-quasiuniform grids (see Figure \ref{fig:counterexample}):


\begin{exa}[Staggered Lax-Friedrichs]
    Consider the trivial problem 
    \[ u_t=0, \qquad u_0=\chi_{[0,1]}. \]
    We discretize it na\"\i vely with the
    staggered Lax-Friedrichs scheme, for flux $f=0$, taking $\Delta_t=h^3$ where $h>0$ is the spatial cell size:
    let $\sys C^h$ contain the cells
    \begin{alignat}{1}
        E^n_j &:= [2nh,(2n+1)h] \times [jh,(j+1)h],\notag\\
        O^n_j &:= [(2n+1)h,(2n+2)h]\times[(j+\ef{2})h,(j+\frac{3}{2})h]\qquad(n\in\Nnull,~j\in\Z). \notag
    \end{alignat}
    Define 
    \begin{alignat}{1}
        F_{E^n_j\rarrow O^n_j}(u^h) &= F_{E^n_j\rarrow O^n_{j-1}} = \frac{h}{2}u^h_{E^n_j},\notag \\
        F_{O^n_j\rarrow E^{n+1}_j}(u^h) &= F_{O^n_j\rarrow E^{n+1}_{j+1}} = \frac{h}{2}u^h_{O^n_j} \qquad(n\in\Nnull,~j\in\Z),\notag \\
        F_{\partial\rarrow E^0_j}(u^h) &= \int_{2jh}^{(2j+1)h}u_0(x)dx \qquad(j\in\Z), \notag
    \end{alignat}
    \begin{alignat}{1}
        F_{E^n_j\rarrow E^n_{j+1}}(u^h) = 0, \qquad F_{O^n_j\rarrow O^n_{j+1}}(u^h) = 0 \qquad(n\in\Nnull,~j\in\Z), \notag \\
        u_{E^0_j} := h^{-1}\int_{jh}{(j+1)h} u_0(x)dS(x). \notag
    \end{alignat}
    It is easy to check that the numerical fluxes are consistent and satisfy the initial condition; all other requirements are satisfied as well. However, 
    the (uniquely determined) solutions $(u^h)_{h>0}$ at a fixed time $t$ approximate $G(\frac{t}{4h},\cdot)*\chi_{[0,1]}$ (as $h\downarrow 0$),
    where $G$ is the heat kernel (this is easy to prove by considering two steps of the numerical scheme: 
    \[ u^h_{E^{n+1}_j} = \frac{u^h_{E^n_{j-1}} + 2u^h_{E^n_j} + u^h_{E^n_{j+1}}}{4} 
    = u^h_{E^n_j} + \frac{\Delta_t}{4h}\frac{u^h_{E^n_{j-1}} - 2u^h_{E^n_j} + u^h_{E^n_{j+1}}}{h^2}; \]
    this is a well-known finite difference scheme for $u_t=\ef{h}u_{xx}$, so standard theory applies).
    However, $\frac{t}{2h}\rightarrow\infty$, so the $u^h$ converge in $L^1_{\loc}$ to $0$, 
    \emph{not} to the actual solution $\chi_{[0,1]}$.
\end{exa}

The example is so ``economical'' that it rules out any conceivable relaxation of the quasiuniformity requirement: the cells are identical rectangles 
(in particular they are convex and have well-behaved surfaces) with sides parallel to the coordinate axes, numerical fluxes through a face depend 
only on one of the directly adjacent cells, the numerical solutions $u^h$ are uniquely defined.
The \emph{only} ``violation'' is that cells become small in time direction faster than in space direction.

The example reflects a problem that does not appear in discretizations of semidiscrete schemes with 
artificial viscosity; while the viscosity coefficient of Lax-Friedrichs type
schemes is $\alpha h^2/\Delta_t$ (for some constant $\alpha$), it is typically $\alpha h$ for semidiscrete methods, so no harm can be done
by choosing $\Delta_t$ too small.

Note that \cite{noelle-mv-irregular} presents a related counterexample, namely the Lax-Friedrichs scheme on a uniform Cartesian (non-staggered) grid
with $\Delta t/h\downarrow 0$ as $h\downarrow 0$. In contrast to staggered Lax-Friedrichs, this counterexample does not serve our purposes: 
for such a scheme, the flux between two intervals in a time step is \emph{not} $O(\Delta t)$, as required by uniform
boundedness \eqref{eq:uniform-boundedness}, but $O(1)$.

\if\doprivate
Such ``anisotropic'' cells that flatten out in some directions are necessary for efficient fine resolution of planar shocks and laminar boundary layers; the
cell count for isotropic cells would grow much faster. However, fine resolution might be too difficult in any case. TODO:
maybe isotropic cells work?
\fi

The restriction to quasiuniform grids is a serious one; results for non-quasiuniform grids are highly desirable because
such grid sequences are produced by adaptive refinement and/or adaptive time integration. However, Theorem \ref{thm:lax-wendroff-ich} 
probably remains true in a special case: tensor products of quasiuniform grids, i.e.
\begin{alignat*}{1}
    \sys C^h &= \bigotimes_{\alpha=1}^s\sys C^h_\alpha := \{ \prod_{\alpha=1}^sC_\alpha : C_\alpha\in\sys C^h \},
\end{alignat*}
where each $\sys C^h_\alpha$ is a grid of the type defined in Section \ref{subsection:grids} (for $\alpha=1$, the grid has to cover $\R^{d_1}_+$, for $\alpha>1$ the grid covers $\R^{d_\alpha}$). The case of semidiscrete schemes can be reduced to the tensor grid case.
These questions will be explored in forthcoming work
(see \cite{elling-lax-wendroff-semidiscrete}).


For arbitrary non-quasiuniform grids, on the other hand, it is necessary to impose stronger conditions on the numerical scheme. 
For example, one could study the error estimators that are used by adaptive schemes to determine where to refine the grid or 
decrease the time step, in order to derive additional smoothness or convergence information about $(u^h)_{h>0}$. 
If sufficiently weak assumptions can be made about a large class of error estimators, it might be possible to derive a 
Lax-Wendroff type result for non-quasiuniform grids that is general enough to be interesting.

\section{Novel applications}

\label{section:novel-applications}

%
%
%

\subsection{Local time stepping}

In order to resolve shocks or contact discontinuities well, it is necessary to refine the grid near them. 
The time step is limited by the CFL condition in small cells near these discontinuities and might be 
unnecessarily small for other parts of the domain. For this reason, it can be efficient to use 
different time steps in different regions; some schemes in this spirit have been proposed in \cite{osher-sanders} or 
\cite{elling-diplom} Chapter 4. Theorem \ref{thm:lax-wendroff-ich} is not limited to spatially unstructured grids; grids can 
be unstructured in space-time, as long as they are quasiuniform in the sense of \eqref{eq:quasi-uniformity}.



\subsection{Moving vertices}

\label{section:moving-vertices}

Another use for the generalized Lax-Wendroff theorem is the large class of numerical methods with unsteady grids (see \cite{elling-diplom}
Section 2.1 for adaptation of classical approximate Riemann solvers to this case). 
These are important because some applications have moving domain boundaries, e.g.\ due to wing flutter or rotating turbine blades. 
Moreover, it is often natural or (for high Mach number supersonic flow) more efficient to use Lagrangian methods (grid vertices 
move along with the fluid). To accomodate these methods is straightforward: instead of a tensor product of time axis partition 
and fixed spatial grid, a grid with moving cells is used; faces are no longer either perpendicular or parallel to the time axes. 
(However, whether Theorem \ref{thm:lax-wendroff-ich} is applicable depends on other details of the scheme as well.)


\subsection{Conservative remapping}

\label{section:conservative-remapping}

In numerical computations, it is sometimes necessary to change grids (\defm{remapping}), for example because the old grid has developed singularities
(especially common for Lagrangian schemes when there is strong vorticity in the flow field). See 
\cite{dukowicz,dukowicz-kodis,grandy,jones-rezoning,dukowicz-baumgardner}
for remapping algorithms and applications.
The remapping step should be conservative, for the same reasons that numerical schemes are conservative, and conservative quantities 
should not be ``transported'' during the remapping step more than necessary. A simple way to achieve this is to set
\[ u_N = V(N)^{-1}\sum_O V(O\cap N)u_O \]
where $N$ is a cell in the new grid, $O$ runs over the cells in the old grid and $u_O,u_N$ are densities of conserved quantities in each. 
(However, to achieve higher orders of accuracy, it might be necessary 
to compute polynomial or spline reconstructions $v(x)$ from the cell averages $u_O$ and to set
\[ u_N := \aint_N v(x)dV(x); \]
the previous scheme corresponds to the obvious piecewise constant reconstruction.)

It is not clear whether  remapping can prevent an otherwise fine numerical scheme from 
converging to the entropy solution. 
However, Theorem \ref{thm:lax-wendroff-ich} can be applied to answer this question for \emph{conservative} remapping: the remapping step 
is interpreted as an extra hyperplane of faces (perpendicular to the time axis), with numerical fluxes defined depending on the remapping algorithm. 
The following requirements are weak enough to cover most existing methods: 
\begin{enumerate}
\item Consistency: whenever all $u_{O'}$ are constant $=w$, the reconstruction $v$ should be constant $=w$ in each cell $O$.
\item Continuity: the reconstruction map $u^h\mapsto v$ should be continuous on the ``diagonal'' of constant grid functions 
    in the $L^\infty\rightarrow L^\infty$ topology.
\item Boundedness: the reconstruction map should be uniformly bounded; more precisely, for any $M$ there should be a constant $c$ so that
    \[ \sup_O|u_O| \leq M \qquad\Rightarrow\qquad \int_O|v(x)|dV(x)\leq cV(O). \]
\item Locality: the values of $v$ over a cell $O$ should depend only on cell averages $u_{O'}$ for $d(O',O)\leq ch$, $c$ some constant independent of $h,O,O'$.
\item Conservation: the reconstruction $v$ should satisfy
    \[ S(O)u^h_O = \int_Ov(x)dV(x). \]
\end{enumerate}
For every remapping step at some time $t$, the old and new cells generate a layer of faces. 
For every cell $N$ that meets a cell $O'$ with $d(O',O)\leq Ch$, add a face $O\rarrow N$ and set 
\[ F_{O\rarrow N} = \aint_N v(x)dV(x). \]

By the assumptions, $F$ has bounded stencil, is consistent, uniformly continuous and uniformly bounded; 
defining $F_{N\rarrow O}:=-F_{O\rarrow N}$ renders it conservative. Moreover, it is clear that if the old and new grid satisfy the requirements 
outlined in the introduction and if the remapping steps are at least $\Omega(h)$ apart, 
then the resulting space-time grid is quasiuniform.

This technique will be discussed in more detail in future work (see \cite{elling-remapping}).

\subsection{Selfsimilar flow}

Selfsimilar solutions, i.e.\ those that satisfy $u(t,x)=u(st,sx)$ for all $s>0$, arise in many important circumstances, such as shocks, contacts or
rarefaction waves Riemann problems, or in \cite{elling-nonuniqueness}. A selfsimilar solution to $u_t+\Div\Vec f(u)=0$
satisfies the system
\begin{alignat*}{1}
    & \Div_\xi(\Vec f(u)-\Vec\xi u) = -du
\end{alignat*}
where $\xi=\frac{x}{t}$ are called \defm{similarity coordinates}. To find asymptotically stable selfsimilar solutions, one can solve
\begin{alignat*}{1}
    & u_\tau + \Div_\xi(\Vec f(u)-\Vec\xi u) = -du
\end{alignat*}
(where $\tau$ is a time-marching ``pseudo-time'' without physical significance). Adding the source term $-du$ to numerical schemes in a conservative
and consistent way (as defined in Section \ref{subsection:numerical-sources}) is easy; achieving stability is not difficult either.
Theorem \ref{thm:lax-wendroff-ich} states that such a scheme delivers entropy solutions as long as it converges.

%
%




\bibliographystyle{amsalpha}
\bibliography{/home/elling/tex/elling}

Preprint Version pre-4, \number\month/\number\day/\number\year

\end{document}